\documentclass{article}
\usepackage{amsmath}
\usepackage{graphicx}
\usepackage{amsfonts}
\usepackage{amssymb}

\begin{document}

\title{Spencer Manifolds}
\author{S. Dimiev$^{+}$, R. Lazov$^{+}$, N. Milev*\\$^{+}$Institute of Mathematics and Informatics,\\Bulg. Acad. of Sciences, \\Acad. G. Bonchev Str., Bl. 8, \\1113 Sofia, Bulgaria\\** {\small Plovdiv University, Math. Faculty, }\\{\small Tzar Assen Str., 24}\\{\small 4000 Plovdiv, \ Bulgaria}}
\date{To appear in ''Quaternionic Structures in Mathematical Physics, Roma, 1999''}
\maketitle
\begin{abstract}
{\small Almost-complex and hyper-complex manifolds are considered in this
paper from the point of view of complex analysis and potential theory. The
idea of holomorphic coordinates on an almost-complex manifold }$(M,\mathbf{J}%
)$ {\small is suggested by D. Spencer [Sp]. For hypercomplex manifolds we
introduce the notion of hyper- holomorphic function and develop some analogous
statements. Elliptic equations are developed in a different way than D.
Spencer . In general here we describe only the formal aspect of the developed theory.}
\end{abstract}

\subsection{INTRODUCTION.}

Differentiable manifolds are described locally by smooth real coordinates.
This is typical in differential geometry. Complex-analytic manifolds are
equipped locally by complex-analytic coordinates. This give rise to the
possibility of applying the theory of holomorphic functions of many complex
variables in the local geometry of complex-analytic manifolds. In the case of
almost complex manifolds $(M,\mathbf{J})$ one use ordinary real coordinates
$(x^{1},...,x^{2n}).$ Here we shall consider complex self-conjugate
coordinates $(z^{1},...,z^{n},\bar{z}^{1},...,\bar{z}^{n})$, where
$z^{k}=x^{2k-1}+ix^{2k}$ and $\bar{z}^{k}=x^{2k-1}-ix^{2k}$.

We denote by $\mathbf{J}^{\ast}$ the action of \textbf{J} on differential
forms of $M$, i. e. by definition

$(\mathbf{J}^{\mathbf{\ast}}(\omega)X\overset{def}{=}\omega(\mathbf{J}X)$,
where $X$ is a vector field, and $\omega$ is a differential form on $M$.

For a fixed index $k$, we say that $z^{k}$ is a ''holomorphic'' coordinate if
$\mathbf{J}^{\ast}dz^{k}=idz^{k}$ and \textbf{J*}$d\bar{z}^{k}=-id\bar{z}^{k}%
$. For non-holomorphic coordinates $z^{q}$ we have
\[
\mathbf{J}^{\mathbf{\ast}}dz^{q}=J_{q}^{1}dz^{1}+....+J_{q}^{n}dz^{n}%
+J_{q}^{n+1}d\bar{z}^{n+1}+....+J_{q}^{2n}d\bar{z}^{2n}%
\]

In the case $z^{k}$ is a holomorphic coordinate for each $k=1,...,n$, the
almost complex structure \textbf{J} is an integrable one. The interest of the
existence of holomorphic coordinates $z^{k}$ when the index $k$ takes not all
values $1,...,n$ is suggested by Donald Spencer [Sp].

By $\mathbb{H=H}(\mathbf{1,}\textsl{i,j,k),ij=k}$ , we will denote the
4-dimensional quaternionic vector space, i. e. $q\in\mathbb{H}$ means that
$q=x^{0}+\textsl{i}x^{1}+\textsl{j}x^{2}+\textsl{k}x^{3},$ where $x^{0}%
,x^{1},x^{2},x^{3}\in\mathbb{R.}$ We will use different complex number
representation for quaternions $q,$ namely $q=z+\zeta\textsl{j,}$ where
$z=x^{0}+\textsl{i}x^{1}$and $\zeta=x^{2}+\textsl{i}x^{3}.$ So we obtain the
right $\textsl{j}$\textsl{-}$\operatorname{complex}$ splitting of
$\mathbb{H},$ denoted by $\mathbb{H}^{j},$ i. e. $\mathbb{H}^{j}%
=(\mathbb{R}\oplus\textsl{i\ }\mathbb{R)\oplus}(\mathbb{R}\oplus
\textsl{i\ }\mathbb{R)}\textsl{j}$ . By $\mathbb{R}\oplus\textsl{i\ }%
\mathbb{R}$ is denoted the tensor product of $\mathbb{R}$ with itself under
the basis ($1,0)$ and $(0,i).$ Identifying $\mathbb{R}\oplus\textsl{i\ }%
\mathbb{R}$ with $\mathbb{C}$ we have that $\mathbb{H}^{j}$ is homeomorphic to
$\mathbb{C\times C}.$ Analogously, we will consider the right \textsl{i-}%
complex splitting of $\mathbb{H},$ namely $\mathbb{H}^{i}=(\mathbb{R}%
\oplus\textsl{j\ }\mathbb{R)\oplus}(\mathbb{R}\ominus\textsl{j\ }%
\mathbb{R)}\textsl{i}$ \ ,\ i.e. $q=x^{0}+\textsl{j}x^{3}+$ +$(x^{2}%
-\textsl{j}x^{4})\textsl{i}$ . $\mathbb{H}^{i}$ is homeomorphic to
$\mathbb{C\times C}$ too.

By $\mathbb{H}^{n}$ is denoted the $n$-dimensional quaternionic vector space
(real $4n$-dimensional)
\[
\mathbb{H}^{n}=\{(q^{1},...,q^{n}):q^{\alpha}\in\mathbb{H},\alpha=1,...,n\}
\]
According to the above accepted notation we have $q^{\alpha}=z^{\alpha}%
+\zeta^{\alpha}\textsl{j}$ , $\alpha=1,...,n$ or
\[
\mathbb{H}^{n}=\mathbb{C}^{n}+\mathbb{C}^{n}\textsl{j}\text{ , \ }%
\mathbb{C}^{n}=\{z^{1},...,z^{n}:z^{\alpha}\in\mathbb{C}\}
\]
This representation is with respect of the right \textsl{j}%
-$\operatorname{complex}$ splitting $\mathbb{H}^{j}.$ A similar representation
of $\mathbb{H}^{n}$ can be written with respect to the right \textsl{i-}%
complex splitting $\mathbb{H}^{i}:\mathbb{H}^{i}=\mathbb{C}^{n}+\mathbb{\bar
{C}}^{n}\textsl{i}$ , \ $\mathbb{C}^{n}=\mathbb{R}^{n}\oplus\mathbb{R}%
^{n}\textsl{i}$ $,$ $\mathbb{\bar{C}}^{n}=\mathbb{R}^{n}\ominus\mathbb{R}%
^{n}\textsl{i}$ .

Let $(M,\mathbf{J},\mathbf{K})$ be a hyper-complex manifold, $\mathbf{JK}%
+\mathbf{KJ}=0,$ $\dim M_{\mathbb{R}}=4n.$ A pair of complex coordinates
$(z,\zeta)$ is called \textit{hyper-holomorphic pair }if\textit{\ }$z$ is
holomorphic with respect to the almost-complex manifold $(M,\mathbf{J})$ and
$\zeta$ is holomorphic with respect to ($M,\mathbf{K}).$

\subsection{Holomorphic coordinates}

\subsubsection{Almost-holomorphic functions}

By definition a function$\ f:U\rightarrow\mathbf{C}$ where $U$ is an open
subset of $M$, is called \textit{almost holomorphic} or almost complex if
$\bar{\partial}f=0$. The above definition can be reformulated in the following
equivalent form:
\[
f\text{\textit{\ is almost holomorphic iff \ }}\mathbf{J}^{\mathbf{\ast}%
}df=idf
\]
\qquad\qquad Respectively, $f$ is \textit{almost-antiholomorphic iff}
$\mathbf{J}^{\mathbf{\ast}}df=-idf$ . For the proof of the equivalence it is
enough to take in view that the exterior derivative $d$ is decomposed as
$d=\partial+\bar{\partial}$ over the space of smooth functions on $M$. Another
form of this definition is obtained taking the real and imaginary parts of
$f$, i.e. $f=u+iv$. In view of $df=du+idv$ we receive $\mathbf{J}%
^{\mathbf{\ast}}du+i\mathbf{J}^{\mathbf{\ast}}dv=idu-dv$. This means that
$\mathbf{J}^{\mathbf{\ast}}du=-dv$ and $\mathbf{J}^{\mathbf{\ast}}dv=du$. As
the obtained two equations are not independent, we can state the following
Cauchy- Riemann type form of the definition

$f=u+iv$ \textit{is almost-holomorphic iff} $\mathbf{J}^{\mathbf{\ast}}dv=du$
\textit{or equivalently} $\mathbf{J}^{\mathbf{\ast}}du=-dv$.

Respectively: $f=u+iv$ \textit{is almost-anti holomorphic iff} $\mathbf{J}%
^{\mathbf{\ast}}dv=-du$ \textit{or equivalently} $\mathbf{J}^{\mathbf{\ast}}du=dv$

\textsl{Remark}$:$ For an almost complex manifold $(M,\mathbf{J})$ with
non-integrable\textbf{\ J}, the decomposition $d=\partial+\bar{\partial}$ is
not valid over differential $(p,q)$-forms on $(M,\mathbf{J})$.

The following proposition is well-known:

\textbf{Proposition 1.} The almost complex structure \textbf{J} of the almost
complex manifold $(M,\mathbf{J})$, $dim_{\mathbf{R}}M=2n$, is an integrable
almost complex structure if and only if for every point $p\in M$, there is a
neighborhood $U$ of $p$ and almost holomorphic functions $f_{j}:U\rightarrow
\mathbb{C},$ $j=1,...,n,$ which differentials at $p$, i.e. $d_{p}f_{j},$
$j=1,...,n,$ are $\mathbb{C}$-linear independent.

\textsl{Remark:} Taking $(U;f_{1},...,f_{n})$ as local coordinate system (as
$f_{j}$ are functionally independent on a neighborhood of $p$), we obtain a
local complex-analytic coordinate system $(U;z_{1},...,z_{n})$, where
$z^{k}=f_{k}$.

\subsubsection{Spencer coordinates}

We say that a local Spencer coordinate system of type $m$ is defined on an
almost complex manifold $(M,\mathbf{J})$ if the following two conditions hold:

1.) There exist an open subset $U$ of $M$ and $m$ different functionally
independent almost holomorphic functions $f_{j}:U\rightarrow\mathbb{C},$
$j=1,...,m,$ such that

2.) The sequence $\ (f_{1},...,f_{m})$ is a maximal sequence of functionally
independent on $U$ almost-holomorphic functions.

3.) The sequence
\[
(U,w^{1},...,w^{m},z^{m+1},...,z^{n},\bar{w}^{n+1},...,\bar{w}^{n+m},\bar
{z}^{n+m+1},...,\bar{z}^{2n})
\]
where $w^{j}=f_{j},$ $j=1,...,m$, determines a local self-conjugate system on
$(M,\mathbf{J})$.

An almost complex manifold which is equipped with an atlas of local Spencer
coordinate systems is by definition an almost-complex manifold of Spencer type
$m$. It is to remark that the notion of Spencer type is correctly defined in
the category of almost complex manifolds. This follows by the fact that each
composition of almost-holomorphic mappings and each inverse of
almost-holomorphic diffeomorphism are almost-holomorphic too.

\textbf{Lemma 1:} The matrix representation of \textbf{J*} in each local
Spencer coordinate system
\[
(U,w^{1},...,w^{m},z^{m+1},...,z^{n},\bar{w}^{n+1},...,\bar{w}^{n+m},\bar
{z}^{n+m+1},...,\bar{z}^{2n})
\]
where $w^{j}=f_{j},$ $j=1,...,m$, are functionally independent almost
holomorphic functions, seems as follows

$\qquad\qquad\qquad\qquad\qquad\qquad\left(
\begin{array}
[c]{cccc}%
iE_{m} & \ast & 0 & \ast\\
0 & \ast & 0 & \ast\\
0 & \ast & -iE_{m} & \ast\\
0 & \ast & 0 & \ast
\end{array}
\right)  $

$E_{m}$ being the unit $m\times m$ matrix.

\textit{Proof. }It is enough to take in view that:
\[
(dw^{1},...,dw^{m},dz^{m+1},...,dz^{n},d\bar{w}^{n+1},...,d\bar{w}^{n+m}%
,d\bar{z}^{n+m+1},...,d\bar{z}^{2n})
\]
is basis of the cotangent space and
\[
\mathbf{J}^{\mathbf{\ast}}dw\/\smallskip^{j}=\mathbf{J}^{\mathbf{\ast}}%
df_{j}=idf_{j}=idw\ ^{j},\text{ }j=1,...,m\text{ }\blacksquare
\]

\textsl{Consequences:} The first $m$ equations of the system $J^{\ast}df=idf$
are just the conditions $\partial f/\partial\bar{z}_{j}=0,$ $j=1,...,m$

We shall consider the mapping from $U$ to $\mathbb{C}^{m}$ defined by
$f_{1},...,f_{m}$. This mapping is a smooth submersion as it can be considered
as a composition of the diffeomorphism defined by Spencer coordinates of $U$
in $\mathbb{C}^{n}\times\mathbb{\bar{C}}^{n}$ and the projection of
$\mathbb{C}^{n}\times\mathbb{\bar{C}}^{n}$ on $\mathbb{C}^{m},$ $m<n$. This
mapping will be denoted by $f_{U}$, and the image of $U$ by $f_{U}$ will be
denoted $U_{m}^{c}.$ It is an open subset of $\mathbb{C}^{m}$, which will be
called a naturally associated $m$-dimensional open set to the considered local
Spencer coordinate system.

\textbf{Lemma 2:} Each almost holomorphic function $h$, defined on a local
Spencer coordinate system $U$ is represented as a superposition of a
holomorphic function $H$ defined on $U_{m}^{c}$ and the almost holomorphic
functions $f_{1},...,f_{m}$ defined on $U$, i.e.
\[
h=H\circ(f_{1},...,f_{m})=H(f_{1},...,f_{m})
\]
\qquad\qquad\qquad

\textit{Proof:} As $w_{j}=f_{j},\ j=1,...,m,$ is a system of smooth
functionally independent on $U$ functions, we have $h=H(w^{1},...,w^{m})$ with
$H\in\mathcal{C}^{\infty}(U)$. But
\[
\bar{\partial}H=(\bar{\partial}H/\partial\bar{w}^{1})d\bar{w}^{1}%
+...+(\partial H/\partial\bar{w}^{m})d\bar{w}^{m}%
\]
and in view of $\bar{\partial}H=\bar{\partial}h=0$, we get that the above
written (0,1)-form is a zero-form, or $\partial H/\partial\bar{w}_{j}$ =0,
$j=1,...m.$ \ $\blacksquare$

\textbf{Lemma 3:} Let $(w^{1},...,w^{m})$ and $(v^{1},...,v^{m})$ be two
systems of holomorphic coordinates on $U_{m}^{c}$ defined by two different
systems of almost holomorphic on $U$ systems $(f_{1},...,f_{m})$ and
$(h_{1},...,h_{m})$. Then there exists a bijective holomorphic transition
mapping between the mentioned two coordinate systems.

\textit{Proof.} According to \textit{Lemma 2} we have $v_{j}=H_{j}%
(w^{1},...,w^{m})$, $j=1,...,n,$where $H_{j}$ are holomorphic functions of
$(w_{1},...,w_{m})$. The system $H=(H_{1},...,H_{m})$ defines the mentioned
transition mapping as the differentials $dH_{j}\ $which are $\mathbb{C}%
$-linear independent. $\blacksquare$

Recapitulating we obtain the following

\textbf{Proposition 2:} On each paracompact almost complex manifold $(M,J)$ of
constant Spencer type $m$ there exists a locally finite covering $U_{j}$ by
self-conjugated Spencer's coordinate system $(Uj,z_{j}^{1},...,z_{j}^{m},...)$
such that in every intersection $U_{j}\cap U_{k}$ the holomorphic coordinates
$z_{j}^{1},...,z_{j}^{m}$ change holomorphically in the other holomorphic
coordinates $z_{k}^{1},...,z_{k}^{m} $.

\subsubsection{Local submersions and local foliations}

As it was remarked above the mapping $f_{U}:U\rightarrow\mathbb{C}^{m}$,
defined by the almost holomorphic functions $(f_{1},...,f_{m})$ is a local
submersion. According to introduced notations
\[
f_{U}(U)=U_{m}^{c}\subset\mathbb{C}^{m}%
\]

The leaves of this submersion are defined as the stalks of the mapping $f_{U}
$. Each leaf is a smooth $(2n-2m)$-dimensional submanifold of $U$ on which all
functions $f_{j}$ have constant value. Transversal leaves are defined as
univalent inverse images of $U_{m}^{c}$, i.e. as sections of $U$ over
$U_{m}^{c}$.

We shall consider the set of all open subsets $U_{m}^{c}$ $\subset
\mathbb{C}^{m}$, corresponding to different mappings $f_{U}$, $U$ open subset
of $M$. This set together with the transition mappings described in
\textit{Lemma 3} defines a pseudo-group of holomorphic transition mappings
between open subsets of $\mathbb{C}^{m}$ denoted as follows
\[
\Gamma\{U_{m}^{c},V_{m}^{c},...;H:U_{m}^{c}\rightarrow V_{m}^{c},....\}
\]

We shall denote by $\mathbf{C}^{m}/\Gamma$ the set of equivalent points of
$\mathbb{C}^{m}$ with respect to the natural equivalence defined by the
holomorphic transition mappings. With this in mind we consider the family
$\{f_{U}:U\rightarrow M\}$ and will define a glued mapping
\[
f:M\rightarrow\mathbb{C}^{m}/\Gamma
\]
as follows: if $p\in M$ we take an open subset $U$ such that $p\in U$ and we
set
\[
f(p)=\{\text{the equivalence class of the point \ }f_{U}(p).\}
\]

Under the assumption that $\mathbb{C}^{m}/\Gamma$ is equipped with the
standard complex structure $\mathbf{i}$ defined by holomorphic coordinates
$(w^{1},...,w^{m})$ we can formulate the following

\textbf{Lemma 4.} The glued mapping $\ f:M\rightarrow\mathbb{C}^{m}/\Gamma$ is
an almost holomorphic mapping between $(M,\mathbf{J})$ and $(\mathbb{C}%
^{m}/\Gamma,\mathbf{i})$.

\textit{Proof.} As the glued mapping $f$ coincides locally with some $f_{U}$
we have:

$\mathbf{J}^{\ast}df_{U}=\mathbf{J}^{\ast}d(f_{1},...,f_{m})=\mathbf{J}^{\ast
}(df_{1},...,df_{m})=(\mathbf{J}^{\ast}df_{1},...,\mathbf{J}^{\ast}%
df_{m})=\mathbf{i}(df_{1},...,df_{m})=id\ f_{U}$. So each $f_{U}$ is an almost
holomorphic mapping $\blacksquare$

\textbf{Lemma 5.} The sheaf of almost holomorphic functions on $M$ is the
inverse image of the sheaf of holomorphic functions on $\mathbb{C}^{m}/\Gamma$.

$Proof.$ The mentioned sheaf on $M$ is defined by the presheaf
$\{U,\mathcal{O}_{M}(U)\}$ where $U$ varies in the set of all open subsets of
$M$ and $\mathcal{O}_{M}(U)$ is defined as follows:

$\mathcal{O}_{M}(U)$ = $\{h\circ f_{U}\mid h\in\mathcal{O}_{\mathbb{C}%
^{m}/\Gamma}f_{U}(U))\}$ .

\subsubsection{Hypercomplex manifolds and hyperholomorphic functions}

Let $M$ be a 4n-dimensional ($\mathcal{C}^{\infty}$) smooth manifold. A
hypercomplex structure on $M$ is defined by a pair of two almost complex
structures \textbf{J} and \textbf{K} such that \textbf{JK + KJ }= 0. It is
easy to see that the composition \textbf{JK} is an almost-complex structure
too. Moreover, for each triple of real numbers $b,c,d$, such that $b^{2}%
+c^{2}+d^{2}=1$, the linear combination $b\mathbf{J}+c\mathbf{K}%
+d(\mathbf{JK})$ is an almost-complex structure on $M$. So there is a family
of almost complex structures on $M$ parametrized by the points of sphere
$\Sigma^{2}$. (See for instance [AM], [ABM]).

We shall consider almost-holomorphic functions on hypercomplex manifolds. The
definition remain the same as in the above considered case, for instance on
$(M,\mathbf{J,K})$ we have\textsl{\ }\textbf{J}-almost- holomorphic function
which are complex-valued function $f$ on $(M,\mathbf{J})$ such that $J^{\ast
}df=idf$ using the right-side \textsl{j}-complex splitting of $\mathbb{H}$.
Respectively \textbf{K}-almost- holomorphic functions $g$ on $(M,\mathbf{J,K}%
)$ are the almost-holomorphic with respect to$(M,\mathbf{K})$ such that
$K^{\ast}dg=jdg$ using an \textsl{i}-complex splitting of $\mathbb{H}$.

Let$(M,\mathbf{J,K})$ be a hypercomplex manifolds and $\mathbb{H}$ be
4-dimensional quaternionic vector space. According to Sommese [So] the
right-side multiplication by $\textsl{i}$ and $\textsl{j}$ are given
respectively by the matrices $S$ and $T$, called standard quaternionic structures$.$

$\ \ \ \ \ \ \ \ \ \ \ \ \ \ S=\left[
\begin{array}
[c]{cccc}%
0 & 1 & 0 & 0\\
-1 & 0 & 0 & 0\\
0 & 0 & 0 & -1\\
0 & 0 & 1 & 0
\end{array}
\right]  $ \ \ \ \ \ \ \ \ \ \ \ \ \ \ $T=\left[
\begin{array}
[c]{cccc}%
0 & 0 & 1 & 0\\
0 & 0 & 0 & 1\\
-1 & 0 & 0 & 0\\
0 & -1 & 0 & 0
\end{array}
\right]  $

In the paper of Sommese the matrix $T$ $is$ denoted by $K.$

As we have $S^{2}=-\mathbf{1},T^{2}=-\mathbf{1},(ST)^{2}=-\mathbf{1}$, and
$ST+TS=0$, we can consider $(\mathbb{H},S,T)$ as a special hypercomplex
manifold. (See [So] ). A function $F$ defined on an open subset $U\subset M$
with valued in $\mathbb{H}$ is called \textbf{J}-hyper-holomorphic function on
$U$ if $dF\circ J=S\circ dF$, or $J^{\ast}dF=SdF$. Using the right-side
\textsl{j}-complex splitting $\mathbb{H}^{j}$ we take the compositions of $F$
with the projections of $\mathbb{H}$ on the first and the second components of
$\mathbb{H}^{j}$. So $F$ is represented by a pair of complex valued functions
denoted respectively by $f$ and $\varphi$. If we set $F=u+\textsl{i\ }v$ +
\textsl{j\ }$\zeta$ + \textsl{k\ }$\eta$, where $u,v,\zeta,\eta$ are
real-valued functions on $U$, we can write $\varphi$ = $u+iv$+ ( $\zeta$
$+i\eta$)\textsl{j}, with $f=u+iv$, $\varphi=\zeta+i\eta$. Complexifying the
matrix $S$, i.e. setting

$\ \ \ \ \ \ \ \ \ \ \ \ \ \ S=$ $\left[
\begin{array}
[c]{cc}%
i & 0\\
0 & -i
\end{array}
\right]  $, \ \ \ \ \ \ \ \ \ \ $\textsl{i}=\left[
\begin{array}
[c]{cc}%
0 & 1\\
-1 & 0
\end{array}
\right]  $, \ \ \ \ \ \ \ \ \ \ \ \ \ $\mathbf{0}=\left[
\begin{array}
[c]{cc}%
0 & 0\\
0 & 0
\end{array}
\right]  $,

and taking $dF=df+d\varphi\textsl{j}$, we calculate that
\[
J^{\ast}df+J^{\ast}d\varphi\textsl{j}=\textsl{i\ }df-\textsl{i\ }%
d\varphi\textsl{j}.
\]
Having in mind the splitting $\mathbb{H}^{j}$, we get $J^{\ast}df=\textsl{i\ }%
df$ and $J^{\ast}d\varphi$ = $-\textsl{i\ }d\varphi$, which means that $f$ is
\textbf{J}-almost-holomorphic function on $U$ and $\varphi$ is \textbf{J}-almost-antiholomorphic.

For the definition of $\mathbf{K}$-hyper-holomorphic function on $U$ we shall
use the other complex splitting of $\mathbb{H}$, namely $\mathbb{H}^{i} $. A
function $G:M\rightarrow\mathbb{H}^{i}$, i.e. $G=g+\psi\textsl{i}$,
$g=u^{\prime}+j\zeta^{\prime}$, $\psi=v^{\prime}-j\eta^{\prime}$, will be
called \textbf{K}-hyper-holomorphic function on $U$ if $\ dG\circ K=T\circ dG
$ or $K^{\ast}dG=TdG$. Taking a (2$\times$2)-representation of the matrix $T$, i.e.

$\ \ \ \ \ \ \ \ T=\left[
\begin{array}
[c]{cc}%
0 & \mathbf{1}\\
-\mathbf{1} & 0
\end{array}
\right]  $, $\ \ \ \ \ \ \ \ \ \mathbf{1}=\left[
\begin{array}
[c]{cc}%
1 & 0\\
0 & 1
\end{array}
\right]  $ , $\ \ \ \ \ \ \ \ 0=$ $\left[
\begin{array}
[c]{cc}%
0 & 0\\
0 & 0
\end{array}
\right]  $, \medskip

after a short calculation we get
\[
K^{\ast}dg+\textsl{i\ }K^{\ast}d\psi=d\psi-\textsl{i\ }dg
\]
It follows that $K^{\ast}dg=d\psi$ and $K^{\ast}d\psi=-dg$. This result is in
terms of $\mathbb{H}^{i}$.

Now we will translate the obtained result in terms of $\mathbb{H}^{j}$. From
$K^{\ast}(du^{\prime}$ + $d\zeta^{\prime}\textsl{j}$) = d$v^{\prime}-$
d$\eta^{\prime}\textsl{j}$ we get
\[
K^{\ast}du^{\prime}=dv^{\prime}\text{ \ and \ \ \ }K^{\ast}d\zeta^{\prime
}=-d\eta^{\prime}.
\]
Analogously, from $K^{\ast}(dv^{\prime}-d\eta^{\prime}\textsl{j}%
)=-(du^{\prime}+d\zeta^{\prime}\textsl{j})$ we get
\[
K^{\ast}dv^{\prime}=-du^{\prime}\text{ \ \ \ and \ \ \ \ }K^{\ast}%
d\eta^{\prime}=-d\zeta^{\prime}.
\]
But the system $K^{\ast}du^{\prime}=dv^{\prime},$ $K^{\ast}dv^{\prime
}=-du^{\prime}$ is just the Cauchy-Riemann system, which says that the
function $u^{\prime}+iv^{\prime}$ is $\mathbf{J}$-almost-antiholomorphic, i.
e. $J^{\ast}d(u^{\prime}+iv^{\prime})=$ =$-id(u^{\prime}+iv^{\prime}).$ The
function $\zeta^{\prime}+i\eta^{\prime}$ is $\mathbf{J}$-almost-holomorphic.

\subsubsection{Hyper-Spencer coordinates}

Hyper-holomorphic coordinates on a hyper-complex manifold $(M,\mathbf{J}%
,\mathbf{K})$ can be introduced by functionally independent
quaternionic-valued functions $f_{\alpha}+\varphi_{\alpha}\textsl{j}$,
$\alpha=$ $1,...,m,$ $m=(1/2)dimM_{\mathbb{R}}$, or by the complex-valued
function ( $f_{\alpha},\varphi_{\alpha}$). We are interested of the
possibility to have $m<(1/2)\dim M_{\mathbb{R}}$. More precisely, a
\textbf{J}-hyper- Spencer coordinate system is defined locally on $M$ as a
maximal system of $m$ functionally independent \textbf{J}-hyper-holomorphic
functions. A hypercomplex manifold equipped with an atlas of local
\textbf{J}-hyper-Spencer coordinate systems is called a hypercomplex manifold
of Spencer type $m$.

Having in mind the interconnection between \textbf{J}-hyper-holomorphic
functions and \textbf{J}-almost-holomorphic ones we derive the analogues of
the \textit{Lemmas 1,2} and \textit{3} of the previous paragraphs. Let us
remark that in view that $f_{\alpha}$ are \textbf{J}-almost- holomorphic, and
$\varphi_{\alpha}$ are \textbf{J}-almost-antiholomorphic, the corresponding
matrix representations of $J^{\ast}$ is as follows (according to \textit{Lemma
1})
\[
\left[
\begin{array}
[c]{cccc}%
iE_{m} & \ast & 0 & \ast\\
0 & \ast & 0 & \ast\\
0 & \ast & -iE_{m} & \ast\\
0 & \ast & 0 & \ast
\end{array}
\right]  \text{ }\times\left[
\begin{array}
[c]{cccc}%
-iE_{m} & \ast & 0 & \ast\\
0 & \ast & 0 & \ast\\
0 & \ast &  iE_{m} & \ast\\
0 & \ast & 0 & \ast
\end{array}
\right]
\]

Analogously, \textbf{K}-hyper-Spencer coordinates can be introduced with the
help of \textbf{K}-hyper holomorphic mappings. The \textit{Proposition 2}
remains valid for \textbf{J-}holomorphic transition functions and
\textbf{K}-holomorphic transition functions. When the transition
transformations are simultaneously \textbf{J- }and \textbf{K-}holomorphic it
follows that they are affine.

Full coordinate systems defined by $m=(1/2)\dim M_{\mathbb{R}}$ functions
which are both \textbf{J} and \textbf{K} hyper-holomorphic lead to
quaternionic manifolds.\bigskip

\subsection{Elliptic Equations}

\subsubsection{Potential structures on almost-complex manifolds}

Let $(M,\mathbf{J})$ be an almost complex manifold. We shall consider the
following globally defined on $M$ Pfaffian form: $\omega=J^{\ast}du$, where
$u=u(p),p\in M$, is a real-valued smooth (at least of class $\mathcal{C}^{2}$)
function. In the case the 1-form $\omega$ is closed, we will say that $\omega$
defines a \textit{potential structure} on the almost complex manifold
$(M,\mathbf{J})$. On each local real coordinate system $(U,x=(x^{k})),x^{k}%
\in\mathbb{R},k=1,...,2n$, we have a matrix representation of \textbf{J}, i.e.
J = $\parallel$ $J_{j}^{k}(x)\parallel$, where $J_{j}^{k}(x)$ are smooth real
functions on $U$. By $J_{j}$ is denoted the $j$-row of the mentioned matrix
and $\nabla u$ is the gradient of $u$. It is easy to see
\[
Jdu=\sum_{q=1}^{2n}(J_{q}\cdot\nabla u)dx^{q}%
\]
where
\[
J_{q}\cdot\nabla u=\sum_{p=1}^{2n}J_{q}^{p}\frac{\partial u}{\partial x^{p}},
\]

For each potential structure on $(M,\mathbf{J})$ the following two statements hold.

\textbf{Consequence 1.} On every simply connected domain $\Omega\subset M$ it
holds that
\[
\int_{\gamma}J^{\ast}du=0
\]
\bigskip for each closed curve $\gamma$ in $\Omega$.

\textbf{Consequence 2.} The following system
\[
\frac{\partial(J_{q}\bullet\nabla u)}{\partial x^{s}}=\frac{\partial
(J_{s}\bullet\nabla u)}{\partial x^{q}},
\]
$s$, $q=1$, $\cdots$, $2n$, is satisfied locally.

\subsubsection{\bigskip Almost pluri-harmonic functions}

By $(M,\mathbf{J},\omega)$ is denoted an almost-complex manifold
$(M,\mathbf{J})$ equipped with potential structure $\omega$. Then the 1-form
$\omega=J^{\ast}du$ is close, and we have $dJ^{\ast}du=0$. In this case we
will say that the function $u$ is an almost-pluriharmonic function. The
interconnection between almost-pluriharmonic functions and almost-holomorphic
ones (with respect to $\mathbf{J}$) is like this one between pluriharmonic
functions and holomorphic ones. This follows directly form the Cauchy-Riemann
equations $J^{\ast}du=-dv,$ $J^{\ast}dv=du$. Clearly the real part $u$ and the
imaginary part $v$ of the almost-holomorphic function $f=u+iv$ are
almost-pluriharmonic functions.

\subsubsection{Elliptic equations on almost-complex manifolds}

We denote by $\triangle_{\mathbf{J}}$ the following differential operator of
second order (in terms of coordinates)
\[
\Delta_{\mathbf{J}}=\sum_{s,p=1}^{2n}A_{sp}\frac{\partial^{2}}{\partial
x^{s}\partial x^{p}}+\sum_{p=1}^{2n}B_{p}\frac{\partial}{\partial x^{p}}%
\]

where
\[
A_{sp}=\sum_{q=1}^{n}(J_{q}^{s}J_{q}^{p}+\delta_{q}^{s}\delta_{q}^{p}),
\]

and
\[
B_{p}=\sum_{s,q=1}^{2n}J_{q}^{s}\left(  \frac{\partial J_{q}^{p}}{\partial
x^{s}}-\frac{\partial J_{s}^{p}}{\partial x^{q}}\right)  ,
\]

$\delta_{q}^{s},\delta_{q}^{p}$ are the Kronecker symbols. Setting
$A_{J}=\Vert A_{sp}\Vert$, we obtain
\[
A_{J}=JJ^{\ast}+E_{2n}%
\]
where $J^{\ast}$ is the transpose of $J$ and $E_{2n}$ is the unity
$2n\times2n$ matrix.

We emphasize here that now we work with real coordinates, but not with complex
self-conjugate ones. However this corresponds to the Spencer type 0. In the
other extreme case of Spencer type $n$ we have complex-analytic (holomorphic)
coordinates. This is the case of complex analytic manifold with the standard
almost-complex structure denoted by $\mathbf{S}^{0}$ (it is different from
\textbf{S }in the previous paragraph)\textbf{. }

\begin{center}
$-\mathbf{S}^{0}=\left[
\begin{array}
[c]{cc}%
0 & 1\\
-1 & 0
\end{array}
\right]  \times...\times\left[
\begin{array}
[c]{cc}%
0 & 1\\
-1 & 0
\end{array}
\right]  $ $\ \ \ $($n$ times)
\end{center}

As $S^{0}(S^{0})^{\ast}=$ $E_{2n}$ \ we get $A_{S^{0}}=2E_{2n}$ and
$\Delta_{\mathbf{S}^{0}}=2\Delta,$ where $\Delta$ is the Laplace operator in
$2n$ real variables.

\textbf{Proposition 3:}$\triangle_{\mathbf{J}}$ is an elliptic differential operator.

\textit{Proof: }It is sufficient to consider the following inequality
\[
\sum_{s=1}^{2n}\sum_{p=1}^{2n}A_{sp}\xi_{s}\xi_{p}=\sum_{q=1}^{2n}\left(
\sum_{s=1}^{2n}J_{q}^{s}\xi_{s}\right)  ^{2}+\sum_{q=1}^{2n}\left(  \sum
_{s=1}^{2n}\delta_{q}^{s}\xi_{s}\right)  ^{2}\geq\sum_{q=1}^{2n}\xi_{q}%
^{2}\text{. \ \ \ \ \ \ }\blacksquare
\]

Considering the PDE
\[
\triangle_{\mathbf{J}}u=0,
\]
we can state the following

\textbf{Theorem :}\textit{Each almost pluriharmonic function }$u$%
\textit{\ satisfies locally the equation }$\triangle_{\mathbf{J}}u=0$

\textit{Proof: }Let $u$ be almost pluriharmonic, i.e. $dJdu=0$, or the 1-form
$J^{\ast}du$ is closed. According to the previous paragraph $u$ satisfies
locally the following system of PDEs
\[
\frac{\partial(J_{q}\bullet\nabla u)}{\partial x^{s}}=\frac{\partial
(J_{s}\bullet\nabla u)}{\partial x^{q}},
\]
$s$, $q=1$, $\cdots$, $2n$. Now replacing
\[
J_{q}\bullet\nabla u=\sum_{p=1}^{2n}J_{k}^{p}\frac{\partial u}{\partial x^{p}%
}\quad\text{ and \ \ \ \ }J_{s}\bullet\nabla u=\sum_{p=1}^{2n}J_{s}^{p}%
\frac{\partial u}{\partial x^{p}}%
\]
in (4) we obtain the system
\[
\sum_{p=1}^{2n}\left(  \frac{\partial\left(  J_{k}^{p}\frac{\partial
u}{\partial x^{p}}\right)  }{\partial x^{s}}-\frac{\partial\left(  J_{s}%
^{p}\frac{\partial u}{\partial x^{p}}\right)  }{\partial x^{k}}\right)  =0,
\]
$k$, $s=1$, $\cdots$, $2n$. Multiplying each of the above written equations by
$J_{q}^{s}$ and summing with respect to $s$ we obtain
\[
\sum_{p=1}^{2n}\sum_{s=1}^{2n}\left(  J_{k}^{p}J_{q}^{s}\frac{\partial^{2}%
u}{\partial x^{s}\partial x^{p}}-J_{q}^{s}J_{s}^{p}\frac{\partial^{2}%
u}{\partial x^{k}\partial x^{p}}\right)  =\sum_{p=1}^{2n}\sum_{s=1}^{2n}%
J_{q}^{s}\left(  \frac{\partial J_{s}^{p}}{\partial x^{k}}-\frac{\partial
J_{k}^{p}}{\partial x^{s}}\right)  \frac{\partial u}{\partial x^{p}}.
\]
As we have
\[
\sum_{s=1}^{2n}J_{q}^{s}J_{s}^{p}=-\delta_{q}^{p}%
\]
and
\[
\frac{\partial^{2}u}{\partial x^{k}\partial x^{p}}=\sum_{s=1}^{2n}\delta
_{k}^{s}\frac{\partial^{2}u}{\partial x^{s}\partial x^{p}},
\]
we obtain
\[
\sum_{p=1}^{2n}\sum_{s=1}^{2n}\left(  J_{k}^{p}J_{q}^{s}+\delta_{q}^{p}%
\delta_{k}^{s}\right)  \frac{\partial^{2}u}{\partial x^{s}\partial x^{p}}%
=\sum_{p=1}^{2n}\sum_{s=1}^{2n}J_{q}^{s}\left(  \frac{\partial J_{s}^{p}%
}{\partial x^{s}}-\frac{\partial J_{k}^{p}}{\partial x^{s}}\right)
\frac{\partial u}{\partial x^{p}}.
\]
Now taking $q=k$ and summing with respect to $k$ we get exactly
\[
\triangle_{\mathbf{J}}u=0.\text{ \ \ }\blacksquare
\]

In the case $\mathbf{J}=\mathbf{S}^{0}$ the above written equation is just the
classical Cauchy-Riemann system.

\textbf{Consequences:}

\textbf{1. }Each almost pluriharmonic function and respectively every almost
holomorphic function of class $\mathcal{C}^{2}$ on a smooth manifold are of
class $\mathcal{C}^{\infty}$\ too.

\textbf{2. }For connected smooth manifolds the maximum principle holds.

\textbf{3. }In the case of real analytic manifold $M,$equipped with
real-analytic structure \textbf{J, }each \textbf{J}-pluriharmonic and each
\textbf{J}-almost-holomorphic function is real analytic.

\textbf{4.} In the case of connected real analytic manifold $M$ with
real-analytic structure \textbf{J }the principle of unicity of the analytic
continuation holds.

\textit{Remark:} This theorem is inspired from the paper [BKW]. The first
announcement is in [DM]

\subsubsection{The equation $dJ^{\ast}d%
\operatorname{u}%
=0$ in terms of vector fields - commutators and anti-commutators}

Applying the well known formula
\[
d\omega(X,Y)=X(\omega(Y))-Y(\omega(X))-\omega([X,Y])\text{, }\omega\text{ is
1-form, }X,Y\text{ are vector fields}%
\]
to the 1-form $\omega=\mathbf{J}du$ we present the equation (2) in terms of
expressions of vector fields, namely
\[
\lbrack X,Y]_{\mathbf{J}}(u)=\mathbf{J}[X,Y](u)
\]
\ \ where $[X,Y]_{\mathbf{J}}\overset{def}{=}X\circ\mathbf{J}Y-Y\circ
\mathbf{J}X.$ It is to remark that $[X,Y]_{\mathbf{J}}$ is not a vector field.
For instance:
\[
\lbrack X,Y]_{\mathbf{J}}(fh)=[X,Y]_{\mathbf{J}}(f)h+f\ [X,Y]_{\mathbf{J}%
}(h)+X(f)(\mathbf{J}Y)(h)-(\mathbf{J}X)(f)Y(h)+X(h)(\mathbf{J}%
Y)(f)-(\mathbf{J}X)(h)Y(f)
\]

\textit{Some properties of }$[X,Y]_{\mathbf{J}}$

Considering the natural splitting
\[
\mathbb{C}TM=T^{1,0}M\oplus T^{0,1}M
\]

we can take the restriction of $[X,Y]_{\mathbf{J}}$ on $T^{1,0}M.$ This means
that
\[
\mathbf{J}X=iX\text{ \ and \ }\mathbf{J}Y=iY
\]
where $X,Y\in T^{1,0}M.$ So we have
\[
\lbrack X,Y]_{\mathbf{J}}=X\circ(iY)-Y\circ(iX)=i[X,Y]
\]
Analogously\ \bigskip%
\[
\lbrack X,Y]_{\mathbf{J}}=(-i)[X,Y]\text{ \ on }T^{0,1}M
\]

Now we take X$\in T^{1,0}M$ and $Y\in T^{0,1}M$
\[
\lbrack X,Y]_{\mathbf{J}}=X\circ(iY)-Y\circ(-iX)=i(X\circ Y+Y\circ
X)=i\{X,Y\}
\]
Here $\{X,Y\}$ denotes the anticommutator of $X$ and $Y.$ Analogously, if
$X\in T^{0,1}M$ and $Y\in T^{1,0}M$:
\[
\lbrack X,Y]_{\mathbf{J}}=-i\{X,Y\}
\]

\subsubsection{Potential structures on hypercomplex manifolds}

On a hypercomplex manifold $(M,\mathbf{J,K})$ we can consider two separate
potential structures, namely
\[
\omega_{1}=\mathbf{J}^{\ast}du\text{ \ \ and \ \ \ }\omega_{2}=\mathbf{K}%
^{\ast}d\zeta
\]
or the sum
\[
\omega=\mathbf{J}^{\ast}du+\mathbf{K}^{\ast}d\zeta
\]

The corresponding almost-pluriharmonic functions $u,v,\zeta,\eta$ satisfy the
equations:
\[
d\mathbf{J}^{\ast}du=d\mathbf{J}^{\ast}dv=0\text{ \ and \ }d\mathbf{J}^{\ast
}d\zeta=d\mathbf{J}^{\ast}d\eta=0
\]

We have also the natural defined elliptic operators $\Delta_{\mathbf{J}}$ and
$\Delta_{\mathbf{K}}.$ According to the proved theorem:
\[
d\mathbf{J}^{\ast}du=d\mathbf{J}^{\ast}dv=0\text{ \ }\Rightarrow\text{
\ }\Delta_{\mathbf{J}}u=\Delta_{\mathbf{J}}v=0
\]

and
\[
d\mathbf{K}^{\ast}d\zeta=d\mathbf{K}^{\ast}d\eta=0\text{ \ }\Rightarrow\text{
\ }\Delta_{\mathbf{K}}\zeta=\Delta_{\mathbf{K}}\eta=0
\]

For the sum $\omega=\mathbf{J}^{\ast}du+\mathbf{K}^{\ast}d\zeta$ a pair of
functions $(u,\zeta)$\ appears, namely the solutions of the following second
order equation:
\[
d\mathbf{J}^{\ast}du+d\mathbf{K}^{\ast}d\zeta=0\text{ }%
\]

In terms of vector fields the above written equations seem as follows
\[
\lbrack X,Y]_{\mathbf{J}}u=\mathbf{J}[X,Y](u)\text{ \ \ \ and \ \ \ }%
[X,Y]_{\mathbf{K}}u=\mathbf{K}[X,Y](u)
\]

\subsection{Generation of almost-complex structures}

\subsubsection{ Remarks on the local equation of almost-holomorphic functions}

Let $(M,\mathbf{J})$ be an almost-complex manifold, $dimM=2n$. Having in mind
the question of the local integration of the equation $\mathbf{J}^{\ast
}df=idf$, we shall examine how ''far away'' a non- integrable almost complex
structure \textbf{J} is from the classical complex structure related with the
standard almost-complex structure \textbf{S}.

Let $p$ be a point of $M$. Taking an open neighborhood $U$ of the point $p$,
small enough, we can accept that $U$ is a neighborhood of the origin in
$\mathbb{R}^{2n}$ ($p$ to be the origin). Now we shall replace\textbf{\ J} by
it matrix representation $J$ on $U$ and $J^{\ast}$ will denote the transposed
matrix. We will use general real coordinates $x=(x^{1},...,x^{2n}%
)\in\mathbb{R}^{2n}$. Let $G$ denotes a non-degenerate $(2n\times2n)$ matrix,
such that $G^{-1}J^{\ast}(0)G=S^{\ast}$, where $S^{\ast}$ is the transposed
matrix of $S$,

$\ \ \ \ \ \ \ \ \ \ \ \ \ \ \ \ \ \ \ \ \ \ S=\left[
\begin{array}
[c]{cc}%
0 & -E_{n}\\
E_{n} & 0
\end{array}
\right]  $ , $E_{n}$ being the unit $n\times n$ matrix.

For $x\in U$ we set:
\[
G^{-1}J(x)G=\left[
\begin{array}
[c]{cc}%
A(x) & B(x)+E_{n}\\
C(x)-E_{n} & D(x)
\end{array}
\right]
\]

$A(x),B(x),C(x),D(x)$ are $n\times n$ matrices.

Clearly we have:

$\left[
\begin{array}
[c]{cc}%
A(x) & B(x)+E_{n}\\
C(x)-E_{n} & D(x)
\end{array}
\right]  =S^{\ast}$ and $A(0)=B(0)=C(0)=D(0)=0\medskip$

Moreover, we have ($G^{-1}J(x)G)^{2}=-E_{2n},$ which implies the following identities:

\qquad$A^{2}(x)+(B(x)+E_{n})(C(x)-E_{n})=-E_{n}$

\qquad$A(x)(B(x)+E_{n})+(B(x)+E_{n})D(x)=0_{n}$

\qquad$(C(x)-E_{n})A(x)+D(x)(C(x)-E_{n})=0_{n}$

\qquad$(C(x)-E_{n})(B(x)+E_{n})+D^{2}(x)=-E_{2n}$

From the last system it follows that locally is valid:

\qquad$A(x)=-(C(x)-E_{n})^{-1}D(x)(C(x)-E_{n})$

\qquad$B(x)+E_{n}=-(C(x)-E_{n})^{-1}(D^{2}(x)+E_{n})$

Indeed, as
\[
\det(C(0)-E_{n})=(-1)^{n}\neq0
\]
the inverse matrix $(C(x)-E_{n})^{-1}$exists in some neighborhood of the
origin \textbf{0}$\in\mathbb{R}^{n}.$

Now lets consider the equation $(J^{\ast}-iE_{2n})df=0.$ It follows that
\[
(G^{-1}J^{\ast}G-iE_{2n})df=0
\]
and also
\[
\left[
\begin{array}
[c]{cc}%
A(x)-iE_{n} & B(x)+E_{n}\\
C(x)-E_{n} & D(x)-iE_{n}%
\end{array}
\right]  df=0
\]

\textsl{Proposition:}The following block matrix identity is valid:\medskip

$\left[
\begin{array}
[c]{cc}%
A(x)-iE_{n} & B(x)+E_{n}%
\end{array}
\right]  =(A(x)-iE_{n})(C(x)-E_{n})^{-1}\left[
\begin{array}
[c]{cc}%
C(x)-E_{n} & D(x)-iE_{n}%
\end{array}
\right]  \medskip$

\textit{Proof:} Let consider the right side of the identity:

$(A(x)-iE_{n})(C(x)-E_{n})^{-1}\left[
\begin{array}
[c]{cc}%
C(x)-E_{n} & D(x)-iE_{n}%
\end{array}
\right]  =$

=$\left[
\begin{array}
[c]{cc}%
A(x)-iE_{n} & (A(x)-iE_{n})(C(x)-E_{n})^{-1}(D(x)-iE_{n})
\end{array}
\right]  $

But:

$(A(x)-iE_{n})(C(x)-E_{n})^{-1}(D(x)-iE_{n})=B(x)+E_{n},$

as $A(x)=-(C(x)-E_{n})^{-1}D(x)(C(x)-E_{n}).$

The last equality becomes:

$(-(C(x)-E_{n})^{-1}D(x)(C(x)-E_{n})-iE_{n})(C(x)-E_{n})^{-1}(D(x)-iE_{n})=$

$=(C(x)-E_{n})^{-1}(-D(x)-iE_{n})(C(x)-E_{n})(C(x)-E_{n})^{-1}(D(x)-iE_{n})=$

= $-(C(x)-E_{n})^{-1}(D(x)+iE_{n})(D(x)-iE_{n})=$

= $-(C(x)-E_{n})^{-1}(D^{2}(x)+iE_{n})=B(x)+E_{n}.$ $\ \blacksquare$

\textsl{Corollary:} The first $n$ equations of the considered system
\[
(J^{\ast}-iE_{2n})df=0
\]
follow from the last $n$ ones. So we obtain that locally this system is
equivalent to the next one:
\[
\left[
\begin{array}
[c]{cc}%
C(x)-E_{n} & D(x)-iE_{n}%
\end{array}
\right]  df=0
\]

or:
\[
\left[
\begin{array}
[c]{cc}%
E_{n} & (C(x)-E_{n})^{-1}(D(x)-iE_{n})
\end{array}
\right]  df=0
\]

Setting $P(x)\overset{def}{=}(C(x)-E_{n})^{-1}D(x)$ and $Q(x)\overset{def}%
{=}(C(x)-E_{n})^{-1},$ we receive the following block matrix form of the
considered equation of almost holomorphic functions:
\[
\left[
\begin{array}
[c]{cc}%
E_{n} & P(x)+iQ(x)
\end{array}
\right]  df=0.
\]

\subsubsection{Local\textbf{\ reconstruction of J by the matrices }%
$P$\textbf{\ and }$Q\bigskip$}

We will use the following equalities:

$C-E_{n}=Q^{-1}$ ; \ $D=Q^{-1}P$ ; \ $A=-QQ^{-1}PQ=-PQ^{-1}$ ;

$B+E_{n}=-Q((Q^{-1}P)^{2}+E_{n})=-PQ^{-1}P-E_{n}.$

The matrix $J$ can be reconstructed as follows:
\[
J=\left[
\begin{array}
[c]{cc}%
-PQ^{-1} & -PQ^{-1}P-Q\\
Q^{-1} & Q^{-1}P
\end{array}
\right]  \text{ \ \ \ \ \ \ \ \ \ \ \ \ \ (*)}%
\]

The mentioned reconstruction (*) can be considered as a generation of the
matrix representation of $J$ on the open set $U$ by the pair of matrices
$(P,Q)$. Denoting by $\mathcal{M}(U,n)$ the algebra of all $(n\times
n)$-matrices equipped with the topology of coordinate convergence, we can
consider the Cartesian product $\mathcal{M}(U,n)\times\mathcal{M}(U,n)$ with
the product topology as a continuous family which generates the set
$\mathcal{J}(U,2n)$ of all $(2n\times2n)$-matrices $J,$ which verify the
matrix equation
\[
J^{2}+E_{2n}=0,
\]
as a kind of moduli space (locally). More precisely, the following proposition holds

\textbf{Proposition 4:}For each $J\in\mathcal{J}(U,2n)$ there is a pair
$(P,Q)\in$ $\mathcal{M}(U,n)\times\mathcal{M}(U,n)$ such that $J$ is generated
by $(P,Q)$ in the sense of the rule (*). Conversely, each pair $(P,Q)$ defines
a $J$ according to the rule (*). Each sequence $(P_{n},Q_{n})$ of elements of
$\mathcal{M}(U,n)\times\mathcal{M}(U,n)$ determines a sequence of elements of
$\mathcal{J}(U,2n)$, and the limit of the second sequence corresponds by the
rule (*) to the limit of the first sequence.

The proof is clear.

\subsubsection{Global reconstruction of $J.$ \bigskip}

The problem of global reconstruction of almost complex structures on a smooth
manifold by an appropriate algebraic objects is much more difficult. It seems
that an approach can be developed on real-analytic almost complex manifold
$(M,J)$ having local matrix representation for $J$ with real-analytic
coefficients. Now we shall consider the sheaf of germs of almost complex
structures, denoted by $\mathcal{J}(M)$, and the sheaf of germs of pairs of
matrices $(P,Q)$. Supposing that each $J$ can be considered as a global
section of the sheaf $\mathcal{J}(M)$, we can develop the rule (*) for germs
of $\mathcal{J}(M)$ and germs of pairs $(P,Q)$ at each point $p\in M$. The set
of global sections of $\mathcal{J}(M)$ must be generated by the sections of
the sheaf of germs of pairs $(P,Q)$.

\textit{Acknowledgment: }The authors are grateful to the organizer of the
\textit{Second Workshop on Quaternionic Structures in Mathematics and Physics}
for the invitation to present this paper.

\bigskip

\begin{center}
REFERENCES
\end{center}

[\textbf{ABM}] Alekssevski D. V., E. Bonan and S. Marchiafava, \textit{On some
structure equations for} \textit{almost-quaternionic hermitian manifolds. //
}K. Sekigawa and S. Dimiev (Eds.),\textit{\ Complex Structures and Vector
Fields}, World Scientific (Singapore), 1995, p. 114 -134

[\textbf{AM}] \ \ Alekssevski D. V. and S. Marchiafava, \textit{Almost
quaternionic hermitian and quasi-k\"{a}hler manifolds.} // K. Sekigawa and S.
Dimiev (Eds.), \textit{Almost Complex Structures}, World Scientific
(Singapore), 1994, p. 150-175

[\textbf{BKW}] Boothby W., S. Kobayashi and H. Wang, \textit{A note on
mappings and automorphisms of almost-complex manifolds.// }Ann. of Math.,
\textbf{77 (}1963), p. 329-334

[\textbf{DM}] \ \ Dimiev S. and O. Muchkarov, \textit{Fonctions
presque-pluriharmoniques. // }C. R. Acad. Bulg. Sci., \textbf{33} (1980), No. 1,

[\textbf{So] \ \ \ \ }Sommese A. J., \textit{Quatenionic manifolds.} // Math.
Ann. \textbf{212} (1975), p. 191- 214

\textbf{[Sp]} \ \ \ Spencer D., \textit{Potential theory and almost complex
manifolds.} // W. Kaplan (Ed). \textit{Lectures on Functions of Complex
Variables}, Univ. of Michigan Press, 1955

Authors addresses:

{\small S. Dimiev\qquad\qquad\qquad\qquad\qquad R. Lazov}

{\small \qquad Institute of Mathematics and Informatic}

{\small \qquad Bulgarian Academy of Sciences}

{\small \qquad8, Acad. G. Boncev Str.}

{\small \qquad1113 Sofia, Bulgaria}

{\small e-mail: sdimiev@math.bas.bg\qquad lazovr@math.bas.bg}

{\small N. Milev}

{\small \qquad Plovdiv University, Math. Faculty}

{\small \qquad24, Tzar Assen Str.}

{\small \qquad4000 Plovdiv, \ Bulgaria }
\end{document}